\documentclass[11pt]{article}

\usepackage{latexsym, graphicx, epsfig, amsmath, amsthm, epsfig, rotating,   amsfonts,amssymb,subfigure,algorithm,algorithmicx,booktabs,bm,color,bigints}
\usepackage{soul}

\newcommand{\pd}[2]{\frac{\partial#1}{\partial#2}}

\begin{document}

\title{Mori-Zwanzig reduced models for uncertainty quantification}
\author{Jing Li, Panos Stinis \\ 
Pacific Northwest National Laboratory \\
Richland, WA 99354} 

\date {}

\maketitle

\begin{abstract}
In many time-dependent problems of practical interest the parameters and/or initial conditions entering the equations describing the evolution of the various quantities exhibit uncertainty. One way to address the problem of how this uncertainty impacts the solution is to expand the solution using polynomial chaos expansions and obtain a system of differential equations for the evolution of the expansion coefficients. We present an application of the Mori-Zwanzig (MZ) formalism to the problem of constructing reduced models of such systems of differential equations. In particular, we construct reduced models for a subset of the polynomial chaos expansion coefficients that are needed for a full description of the uncertainty caused by uncertain parameters or initial conditions. 

Even though the MZ formalism is exact, its straightforward application to the problem of constructing reduced models for estimating uncertainty involves the computation of memory terms whose cost can become prohibitively expensive. For those cases, we present a Markovian reformulation of the MZ formalism which can lead to approximations that can alleviate some of the computational expense while retaining an accuracy advantage over reduced models that discard the memory altogether. Our results support the conclusion that successful reduced models need to include memory effects. 

\end{abstract}

\section{Introduction}
The problem of quantifying the uncertainty of the solution of systems of partial or ordinary differential equations has become in recent years a rather active area of research. The realization that more often than not, for problems of practical interest, one is not able to determine the parameters, initial conditions, boundary conditions etc. to within high enough accuracy, has led to a flourishing literature of methods for quantifying the impact that this uncertainty imposes on the solution of the problems under investigation (see e.g. \cite{ghanem,leonenko,ma,nouy,venturi,wan,barajas}). However, despite the increase in computational power and the development of various techniques for uncertainty quantification there is still a wealth of problems where reliable uncertainty quantification is beyond reach. One way to address this problem is to look for reduced models for a subset of the variables needed for a complete description of the uncertainty. The effect of all types of uncertainty is intimately connected with the inherent instabilities that may be present in the underlying system which we subject to the uncertainty. These considerations remain equally, if not more, important when we attempt to construct reduced models for uncertainty quantification. 

In the current work, we are concerned with the construction of reduced models for systems of differential equations that arise from polynomial chaos expansions of solutions of a PDE or ODE system. In particular, we focus on the case that the given PDE or ODE system contains uncertain parameters or initial conditions and we want to construct a reduced model for the evolution of a subset of the polynomial chaos expansions that are needed for a complete description of the uncertainty caused by the uncertain parameters. There are different methods to construct reduced models for PDE or ODE systems (see e.g. \cite{givon,CS05} and references therein). We choose to use the Mori-Zwanzig formalism in order to construct the reduced model \cite{CHK00,CHK3}. 

The main issue with all model reduction approaches is the computation of the memory caused by the process of eliminating variables from the given system (referred to as the full system from this point on) \cite{CS05}. The memory terms are, in general, integral terms which account for the history of the variables that are not resolved. In principle the integrands appearing in the memory terms can be computed through the solution of the orthogonal dynamics equation \cite{CHK3}. We present some examples where this procedure can be implemented and the resulting reduced model can be estimated. Those examples highlight the definite improvement in accuracy of a reduced model when it includes a memory term. However, it is also easy to come up with examples where the solution of the orthogonal dynamics equation becomes prohibitively expensive.

For such cases we present a Markovian reformulation of the MZ formalism which allows the calculation of the memory terms through the solution of ordinary differential equations instead of the computation of convolution integrals as they appear in the original formulation. We present an algorithm which allows the estimation of the necessary parameters on the fly. This means that one starts evolving the full system and use it to estimate the reduced model parameters. Once this is achieved, the simulation continues by evolving {\it only} the reduced model with the necessary parameters set equal to their estimated values from the first part of the algorithm. Of course, such an approximation of the memory term cannot work under all circumstances. We present results for a nontrivial problem where it does yield a reduced model with improved behavior compared to a model that ignores the memory terms altogether. 

We should note that this alternative approach to computing the memory term fits in the renormalization framework advocated recently by one of the authors \cite{s11} in order to construct reduced models for singular PDEs. In particular, the idea is that one embeds the MZ reduced model in a larger family of reduced models which share the same functional form but may have additional parameters for enhanced flexibility. These extra parameters are determined so that the reduced model reproduces some dynamic features of the full system. After this is done, one can switch to the reduced model for the rest of the simulation. In the current work, the extra parameters are the lengths of the memory appearing in the MZ reduced model. 

Section \ref{mz_formalism} presents a brief introduction to the MZ formalism for the construction of reduced models of systems of ODEs. In Section \ref{reformulation} we develop the Markovian reformulation of the MZ formalism and show how one can estimate adaptively the parameters appearing in the reduced model. Section \ref{examples} presents numerical results both for the original MZ formalism (Sections \ref{example1}-\ref{example2}) and its Markovian reformulation (Section \ref{example3}). Finally, in Section \ref{discussion} we discuss directions for future work.


\section{Mori-Zwanzig formalism}\label{mz_formalism}

We begin with a brief presentation of the Mori-Zwanzig formalism \cite{CHK00,CHK3}. Suppose we are given the system 
\begin{equation}\label{odes}
\frac{du(t)}{dt} = R (t,u(t)),
\end{equation}
where $u = ( \{u_k\}), \; k \in H \cup G$ 
with initial condition $u(0)=u_0.$ The unknown variables (modes) are divided into two groups, one group is indexed in H and the order indexed in G. Our goal is to construct a reduced model for the modes in the set $H.$ The system of ordinary differential equations
we are given can be transformed into a system of  linear
partial differential equations
\begin{equation}
\label{pde}
\pd{\varphi_k}{t}=L \varphi_k, \qquad \varphi_k (u_0,0)=u_{0k}, \, k \in H \cup G
\end{equation}
where $L=\sum_{k \in H \cup G } R_i(u_0) \frac{\partial}{\partial u_{0i}}.$ The solution of \eqref{pde} is
given by $u_k (u_0,t)=\varphi_k(u_0,t)$. Using semigroup notation we can rewrite (\ref{pde}) as
$$\pd{}{t} e^{tL} u_{0k}=L e^{tL} u_{0k}$$
Suppose that the vector of initial conditions can be divided as $u_0=(\hat{u}_0,\tilde{u}_0),$ where 
$\hat{u}_0$ is the vector of the resolved variables (those in $H$) and $\tilde{u}_0$ is the vector of the unresolved variables (those in $G$).  Let $P$ be an orthogonal projection on the space of functions of $\hat{u}_0$ and $Q=I-P.$ 

Equation \eqref{pde} 
can be rewritten as 
\begin{equation}
\label{mz}
\frac{\partial}{\partial{t}} e^{tL}u_{0k}=
e^{tL}PLu_{0k}+e^{tQL}QLu_{0k}+
\int_0^t e^{(t-s)L}PLe^{sQL}QLu_{0k}ds, \, k \in H,
\end{equation}
where we have used Dyson's formula
\begin{equation}
\label{dyson1}
e^{tL}=e^{tQL}+\int_0^t e^{(t-s)L}PLe^{sQL}ds.
\end{equation}
Equation (\ref{mz}) is the Mori-Zwanzig identity. 
Note that
this relation is exact and is an alternative way
of writing the original PDE. It is the starting
point of our approximations. Of course, we
have one such equation for each of the resolved
variables $u_k, k \in H$. The first term in (\ref{mz}) is
usually called Markovian since it depends only on the values of the variables
at the current instant, the second is called ``noise" and the third ``memory". 

If we write
$$e^{tQL}QLu_{0k}=w_k,$$ 
$w_k(u_0,t)$ satisfies the equation
\begin{equation}
\label{ortho}
\begin{cases}
&\frac{\partial}{\partial{t}}w_k(u_0,t)=QLw_k(u_0,t) \\ 
& w_k(u_0,0) = QLx_k=R_k(u_0)-(PR_k)(\hat{u_0}). 
\end{cases} 
\end{equation}
If we project (\ref{ortho}) we get
$$P\frac{\partial}{\partial{t}}w_k(u_0,t)=
PQLw_k(u_0,t)=0,$$
since $PQ=0$. Also for the initial condition
$$Pw_k(u_0,0)=PQLu_{0k}=0$$
by the same argument. Thus, the solution
of (\ref{ortho}) is at all times orthogonal
to the range of $P.$ We call
(\ref{ortho}) the orthogonal dynamics equation. Since the solutions of 
the orthogonal dynamics equation remain orthogonal to the range of $P$, 
we can project the Mori-Zwanzig equation (\ref{mz}) and find
\begin{equation}
\label{mzp}
\frac{\partial}{\partial{t}} Pe^{tL}u_{0k}=
Pe^{tL}PLu_{0k}+
P\int_0^t e^{(t-s)L}PLe^{sQL}QLu_{0k} ds.
\end{equation}
We will not present here more details about how to start from Eq. \eqref{mzp} and construct reduced models of different orders for a general system of ODEs. Such constructions have been documented thoroughly elsewhere (see e.g. \cite{CHK3}). However, we will provide such details for the specific numerical examples in Sections \ref{example1}-\ref{example2}. 


\section{Markovian reformulation of the MZ formalism}\label{reformulation}
While the MZ model given by Eq. \eqref{mzp} is exact, its construction can be involved and most importantly, very costly. The main source of computational expense is the memory term. Technically, the cost associated with the memory term comes from two sources: i) the presence of the orthogonal dynamics equation solution operator $e^{sQL}$ and ii) the need to find an expression in terms of the resolved variables and time for $PLe^{sQL}QLu_{0k}$ which appears in the memory integrand. 
The presence of $e^{sQL}$ is problematic because the orthogonal dynamics equation is, for the general case, a PDE in as many dimensions as the original system of ODEs. Also, finding an expression for $PLe^{sQL}QLu_{0k}$ is problematic because, in general, it is not possible to separate the dependence of the expression on time and on the resolved variables. Both are formidable tasks and we will show with several examples how they can increase the cost of constructing the reduced model. For some cases (see e.g. Section \ref{example1} and \ref{example2}) both tasks can be tackled through the use of a finite-rank projection for the operator $P.$ However, we will show with a simple example (see Section \ref{example3}) that the use of a finite-rank projection may be too costly itself. For such cases, we need an alternative approach to the construction of the memory term. In this section we describe a reformulation of the problem of computing the memory term which 
can alleviate some of these issues. Also, we present numerical results from the application of this approach in Section \ref{example3}.

\subsection{Finite memory}\label{memory_comp}

We focus on the case when the memory has a finite extent only. The case of infinite memory is simpler and is a special case of the formulation presented below. Also, the current reformulation allows us to comment on what happens in the case when the memory is very short. 

Let $w_{0k}(t)=P\int_0^t e^{(t-s)L}PLe^{sQL}QLu_{0k} ds=P\int_0^t e^{sL}PLe^{(t-s)QL}QLu_{0k} ds,$ by the change of variables $t'=t-s.$ Note, that $w_{0k}$ depends both on $t$ and the resolved part of the initial conditions $\hat{u}_0.$ We have suppressed the $\hat{u}_0$ dependence for simplicity of notation. If the memory extends only for $t_0$ units in the past (with $t_0 \leq t,$) then $$w_{0k}(t)=P\int_{t-t_0}^t e^{sL}PLe^{(t-s)QL}QLu_{0k} ds.$$ The evolution of 
$w_{0k}$ is given by 
\begin{equation}\label{memory_1}
\frac{dw_{0k}}{dt}=Pe^{tL}PLQLu_{0k}-Pe^{(t-t_0)L}PLe^{t_0 QL}QLu_{0k}+w_{1k}(t),
\end{equation}
where $$w_{1k}(t)=P\int_{t-t_0}^t e^{sL}PLe^{(t-s)QL}QLQLu_{0k} ds.$$ To allow for more flexibility, let us assume that the integrand in the formula for $w_{1k}(t)$ contributes only for $t_1$ units with $t_1 \leq t_0.$ Then $$w_{1k}(t)=P\int_{t-t_1}^t e^{sL}PLe^{(t-s)QL}QLQLu_{0k} ds.$$ 
We can proceed and write an equation for the evolution of $w_{1k}(t)$ which reads
\begin{equation}\label{memory_2}
\frac{dw_{1k}}{dt}=Pe^{tL}PLQLQLu_{0k}-Pe^{(t-t_1)L}PLe^{t_1 QL}QLQLu_{0k}+w_{2k}(t),
\end{equation}
where $$w_{2k}(t)=P\int_{t-t_1}^t e^{sL}PLe^{(t-s)QL}QLQLQLu_{0k} ds.$$ Similarly, if this integral extends only for $t_2$ units in the past with $t_2 \leq t_1,$ then
$$w_{2k}(t)=P\int_{t-t_2}^t e^{sL}PLe^{(t-s)QL}QLQLQLu_{0k} ds.$$
This hierarchy of equations continues indefinitely. Also, we can assume for more flexibility that at every level of the hierarchy we allow the interval of integration for the integral term to extend to fewer or the same units of time than the integral in the previous level. If we keep, say, $n$ terms in this hierarchy, the equation for $w_{(n-1)k}(t)$ will read 
\begin{gather}\label{memory_n}
\frac{dw_{(n-1)k}}{dt}=Pe^{tL}PL(QL)^{n-1}QLu_{0k}- \\
 Pe^{(t-t_{n-1})L}PLe^{t_{n-1} QL}(QL)^{n-1}QLu_{0k}+w_{nk}(t)  \notag
\end{gather}
where $$w_{nk}(t)=P\int_{t-t_n}^t e^{sL}PLe^{(t-s)QL}(QL)^{n}QLu_{0k} ds$$
Note that the last term in \eqref{memory_n} involves the unknown evolution operator for the orthogonal dynamics equation. This situation is the well-known closure problem. We can stop the hierarchy at the $n$th term by assuming that $w_{nk}(t)=0.$

In addition to the closure problem, the unknown evolution operator for the orthogonal dynamics equation appears in the equations for the evolution of the quantities $w_{0k}(t),$ $\ldots,$ $w_{(n-1)k}(t)$ through the various terms $Pe^{(t-t_0)L}PLe^{t_0 QL}QLu_{0k},$ $\ldots,$ $Pe^{(t-t_0)L}PLe^{t_0 QL}(QL)^{n-1}QLu_{0k}$ respectively.

We describe now a way to express these terms involving the unknown orthogonal dynamics operator through known quantities so that we obtain a closed system for the evolution of $w_{0k}(t),\ldots,w_{(n-1)k}(t).$

Since we want to treat the case where $t_0$ is not necessarily small, we divide the interval $[t-t_0,t]$ in $n_0$ subintervals. Define 

\begin{align*}
w_{0k}^{(1)}(t) & =P\int_{t-\Delta t_0}^t e^{sL}PLe^{(t-s)QL}QLu_{0k} ds \\
w_{0k}^{(2)}(t) & =P\int_{t-2 \Delta t_0}^{t- \Delta t_0} e^{sL}PLe^{(t-s)QL}QLu_{0k} ds \\
\ldots & \\
w_{0k}^{(n_0)}(t) & =P\int_{t-t_0}^{t- (n_0-1)\Delta t_0} e^{sL}PLe^{(t-s)QL}QLu_{0k} ds,
\end{align*}
where $n_0 \Delta t_0 = t_0$ and $w_{0k}(t)=\sum_{i=1}^{n_0} w_{0k}^{(i)}(t).$ Similarly, we can define the quantities $w_{1k}^{(1)}(t),\ldots,w_{1k}^{(n_1)}(t)$ 
\begin{align*}
w_{1k}^{(1)}(t) & =P\int_{t-\Delta t_1}^t e^{sL}PLe^{(t-s)QL}QLQLu_{0k} ds \\
w_{1k}^{(2)}(t) & =P\int_{t-2 \Delta t_1}^{t- \Delta t_1} e^{sL}PLe^{(t-s)QL}QLQLu_{0k} ds \\
\ldots & \\
w_{1k}^{(n_1)}(t) & =P\int_{t-t_1}^{t- (n_1-1)\Delta t_1} e^{sL}PLe^{(t-s)QL}QLQLu_{0k} ds,
\end{align*}
where $n_1 \Delta t_1 = t_1$ and $w_{1k}(t)=\sum_{i=1}^{n_1} w_{1k}^{(i)}(t).$ In a similar fashion we can define corresponding quantities for all the memory terms up to  $w_{(n-1)k}(t)=\sum_{i=1}^{n_{n-1}} w_{(n-1)k}^{(i)}(t).$ 

In order to proceed we need to make an approximation for the integrals over the subintervals.

\subsection{Trapezoidal rule approximation}\label{trapezoidal}
We have
\begin{multline*}
w_{0k}^{(1)}(t)  =P\int_{t-\Delta t_0}^t e^{sL}PLe^{(t-s)QL}QLu_{0k} ds   \\
=\biggl[ Pe^{tL}PLQLu_{0k}+Pe^{(t-\Delta t_0)L}PLe^{\Delta t_0 QL}QLu_{0k} \biggr] \frac{\Delta t_0}{2}+ O((\Delta t_0)^3)
\end{multline*}
from which we find
$$Pe^{(t-\Delta t_0)L}PLe^{\Delta t_0 QL}QLu_{0k}=\biggl ( \frac{2}{\Delta t_0} \biggr ) w_{0k}^{(1)}(t) - Pe^{tL}PLQLu_{0k} + O((\Delta t_0)^2)$$
and from \eqref{memory_1}
\begin{equation*}
\frac{dw_{0k}^{(1)}}{dt}=-\biggl ( \frac{2}{\Delta t_0} \biggr ) w_{0k}^{(1)}(t)+ 2Pe^{tL}PLQLu_{0k}+w_{1k}^{(1)}(t)+ O((\Delta t_0)^2).
\end{equation*}
Similarly, for $w_{0k}^{(2)}(t)$ we find
\begin{multline*}
\frac{dw_{0k}^{(2)}}{dt}=\biggl ( \frac{4}{\Delta t_0} \biggr ) w_{0k}^{(1)}(t) \\
-\biggl ( \frac{2}{\Delta t_0} \biggr ) w_{0k}^{(2)}(t) - 2Pe^{tL}PLQLu_{0k} 
+w_{1k}^{(2)}(t)+ O((\Delta t_0)^2)
\end{multline*}
In general,
\begin{multline}\label{memory_1a}
\frac{dw_{0k}^{(i)}}{dt}= -\biggl ( \frac{2}{\Delta t_0} \biggr ) w_{0k}^{(i)}(t) + (-1)^{i+1} 2Pe^{tL}PLQLu_{0k} \\
 +\biggl [  \sum_{j=1}^{i-1}  \biggl ( \frac{4}{\Delta t_0} \biggr ) (-1)^{i+j+1} w_{0k}^{(j)}(t) \biggr ] +w_{1k}^{(i)}(t)+ O((\Delta t_0)^2)  \; \; \text{for}  \; \; i=1,\ldots,n_0.
\end{multline}
Similarly,
\begin{multline*}
\frac{dw_{1k}^{(i)}}{dt}= -\biggl ( \frac{2}{\Delta t_1} \biggr ) w_{1k}^{(i)}(t) + (-1)^{i+1} 2Pe^{tL}PLQLQLu_{0k} \\
 +\biggl [  \sum_{j=1}^{i-1}  \biggl ( \frac{4}{\Delta t_1} \biggr ) (-1)^{i+j+1} w_{1k}^{(j)}(t) \biggr ] +w_{2k}^{(i)}(t)+ O((\Delta t_1)^2)  \; \; \text{for}  \; \; i=1,\ldots,n_1 
\end{multline*} 
$\ldots$
\begin{multline}
\frac{dw_{(n-1)k}^{(i)}}{dt}= -\biggl ( \frac{2}{\Delta t_{n-1}} \biggr ) w_{(n-1)k}^{(i)}(t) + (-1)^{i+1} 2Pe^{tL}PL(QL)^{n-1}QLu_{0k} \\
 +\biggl [  \sum_{j=1}^{i-1}  \biggl ( \frac{4}{\Delta t_{n-1}} \biggr ) (-1)^{i+j+1} w_{(n-1)k}^{(j)}(t) \biggr ] + O((\Delta t_{n-1})^2)  \; \; \text{for}  \; \; i=1,\ldots,n_{n-1}.
\end{multline}
By dropping the $O((\Delta t_0)^2),\ldots, O((\Delta t_{n-1})^2)$ terms we obtain a system of $n_0+n_1+\ldots+n_{n-1}$ differential equations for the evolution of the quantities $w_{0k}^{(1)}(t),\ldots,w_{(n-1)k}^{(n_{n-1})}.$ This system allows us to determine the memory term $w_{0k}(t).$ Since the approximation we have used for the integral leads to an error $O(\Delta t)^2,$ the ODE solver should also be $O(\Delta t)^2.$ We have used the modified Euler method to solve numerically the equations for the reduced model. 

Note that the implementation of the above scheme requires the knowledge of the expressions for $Pe^{tL}PLQLu_{0k},\ldots,Pe^{tL}PL(QL)^{n-1}QLu_{0k}.$ Since the computation of these expressions for large $n$ can be rather involved for nonlinear systems (see Section \ref{example3}), we expect that the above scheme will be used with a small to moderate value of $n.$ Finally, we mention that the above construction can be carried out for integration rules of higher order e.g. Simpson's rule.

\subsection{Estimation of the memory length}\label{mz_length}

The construction presented above relies on an accurate determination of the memory lengths $ t_0,  t_1,\ldots, t_{n-1}.$ We present in this section a way to estimate these quantities on the fly. This means that we start evolving the {\it full} system, use it to estimate $ t_0,  t_1,\ldots, t_{n-1}$ and then switch to the reduced model with the estimated values for $ t_0,  t_1,\ldots, t_{n-1}.$

For simplicity of presentation we assume that we evolve only $w_{0k}(t).$ If we use the trapezoidal rule to discretize $w_{0k}(t)$ and eliminate the term $Pe^{(t-t_0)L}PLe^{t_0 QL}QLu_{0k}$ from \eqref{memory_1}, the reduced model reads 

\begin{gather}
\frac{d Pu_{k}}{dt}= Pe^{tL}PLu_{0k} + w_{0k}(t) \label{reduced1} \\
\frac{d w_{0k}}{dt}=2Pe^{tL}PLQLu_{0k} - \frac{2}{t_0} w_{0k}(t)  \label{reduced2}
\end{gather}
for $k \in H. $ 
We can solve \eqref{reduced2} formally and substitute in \eqref{reduced1} to get
\begin{equation}\label{reduced_integral}
\frac{d Pu_{k}}{dt}= Pe^{tL}PLu_{0k} + \int_0^t e^{-\lambda_0(t-s)}2Pe^{sL}PLQLu_{0k} ds 
\end{equation}
where $\lambda_0=2/t_0.$ Recall that, for the resolved variables, we have from the full system 
\begin{equation}\label{full_split}
\frac{d Pu_{k}}{dt}= Pe^{tL}PLu_{0k} + Pe^{tL}QLu_{0k}.  
\end{equation}
We would like to estimate the memory decay parameter $t_0$ so that the reduced equation \eqref{reduced_integral} for $u_{k}$ reproduces the behavior of $u_{k}$ as predicted by the full system \eqref{full_split}. We can do that by requiring that the evolution of some integral quantity of the solution is the same when predicted by the reduced and full systems. 

We begin by discretizing the integral term in \eqref{reduced_integral}. Suppose that we are evolving the full system with a step size $\delta t,$ where $t=n_t \delta t$ (note that $n_t$ increases as $t$ increases). If we discretize the integral with the trapezoidal rule we find
\begin{gather}\label{reduced_integral2}
\frac{d Pu_{k}}{dt}= Pe^{tL}PLu_{0k}  \\
+ [f_{k}(t,\hat{u}_{0})+2 \sum_{j=1}^{n_t-1}e^{-\lambda_0(t-j\delta t)}f_{k}(j\delta t, \hat{u}_{0}) +e^{-\lambda_0t}f_{k}(0,\hat{u}_{0})] \frac{\delta t}{2} \notag
\end{gather} 
where $f_{k}(j\delta t, \hat{u}_{0})=2Pe^{j\delta t L}PLQLu_{0k}$ for $j=0,\ldots,n_t.$ The quantities $f_{k}(j\delta t, \hat{u}_{0})$ can be computed from the full system.  

There is freedom in the choice of the integral quantity whose evolution the reduced model should be able to reproduce. For example, we can use $ \sum_{k \in H}  |Pu_{k}(t)|^2$ the squared $l_2$ norm of the resolved variables. If we use this integral quantity, then from \eqref{reduced_integral2} and \eqref{full_split} we find that the unknown parameter $t_0$ must satisfy 
\begin{equation}\label{newton1}
 \sum_{k \in H} 2 Re \{  I_{k}(t,t_0) (Pu_{k})^*(t) \} =   \sum_{k \in H}  2 Re \{  Pe^{tL}QLu_{0k} (Pu_{k})^*(t) \} , 
\end{equation} 
where $$I_{k}(t,t_0)= [f_{k}(t,\hat{u}_{0})+2 \sum_{j=1}^{n_t-1}e^{-\lambda_0(t-j\delta t)}f_{k}(j\delta t, \hat{u}_{0}) +e^{-\lambda_0t}f_{k}(0,\hat{u}_{0})] \frac{\delta t}{2}$$ and $Re\{\cdot\}$ denotes the real part.

Let $y=\exp[-\lambda_0\delta t].$ Then, 
\begin{equation}\label{newton2}
I_{k}(t,t_0)= [f_{k}(t,\hat{u}_{0})+2 \sum_{j=1}^{n_t-1} y^{n_t-j}  f_{k}(j\delta t, \hat{u}_{0}) +y^{n_t}f_{k}(0,\hat{u}_{0})] \frac{\delta t}{2}.
\end{equation} 
With this identification, equation \eqref{newton1} becomes a polynomial equation for $y$ with $y \in [0,1].$ It is not difficult to solve equation \eqref{newton1} with an iterative method, for example Newton's method. For the numerical results we present in Section \ref{example3}, Newton's method converged to double precision accuracy within 4-5 iterations. After an estimate $\hat{y}$ has been obtained, we can find the estimate $\hat{t}_0$ of $t_0$ (recall $\lambda_0=2/t_0$) from 
\begin{equation}\label{newton3}
\hat{t}_0=-\frac{2 \delta t}{\ln \hat{y}} .
\end{equation}

\subsubsection{Determination of optimal estimate $\hat{t}_0$}\label{mz_optimal}
For each time instant $t$ we can obtain through equations \eqref{newton1} and \eqref{newton3}, an estimate $\hat{t}_0(t)$ for $t_0.$ Thus, the most important issue that we have to address is that of deciding which is the best estimate of $t_0.$ In other words, at what time $t_f$ should we stop estimating the value of $t_0$ so that we can use the estimated value $\hat{t}_0(t_f)$ to evolve the reduced model from then on. 

We define $\epsilon(t)=\underset{l\in [1,n_t]}{\max} |\hat{y}^l (t+\delta t)-\hat{y}^l (t)|.$ The quantity $\epsilon(t)$ monitors the convergence of not only the value of the estimate $\hat{y}$ as a function of the time $t$, but of the whole function $e^{-\lambda_0(t-s)}.$ Ideally, $\epsilon(t)$ converges to zero with increasing $t.$ That will be the case if the approximation of the memory term only through $Pe^{tL}PLQLu_{0kr}$ is enough (see \eqref{reduced1}-\eqref{reduced2}). However, this will not always be the case. If keeping $Pe^{tL}PLQLu_{0kr}$ is not enough, then $\epsilon(t)$ will decrease with increasing $t$ up to some time $t_{min}$ when it will reach a nonzero minimum. After that time, it starts increasing. This signals that keeping only $Pe^{tL}PLQLu_{0kr}$ is {\it not enough} to describe accurately the memory. 

In order to proceed we have two options: (i) construct a higher order model and (ii) identify $t_f=t_{min}$ and thus $\hat{t}_0(t_f)=\hat{t}_0(t_{min}).$ Results for higher order models will be presented elsewhere (see also discussion in Section \ref{discussion}). In the numerical experiments we present in the next section we have chosen $\hat{t}_0(t_f)=\hat{t}_0(t_{min}).$ Note that the procedure just outlined allows the automation of the algorithm. This means that there is no adjustable reduced model parameter that needs to be specified at the onset of the algorithm. 

We are now in a position to state the adaptive Mori-Zwanzig algorithm which constructs a reduced model with the necessary memory term parameter $t_0$ estimated on the fly. 

\vskip14pt
{\bf Adaptive Mori-Zwanzig Algorithm}
\begin{enumerate}
 \item
Evolve the full system and compute, at every step, the estimate $\hat{t}_0(t).$ Use estimates of $t_0$ from successive steps to calculate $\epsilon(t)=\underset{l\in [1,n_t]}{\max} |\hat{y}^l (t+\delta t)-\hat{y}^l (t)|.$  
\item
When $\epsilon(t)$ reaches a minimum (possibly non zero) value at some instant $t_{min}$, pick $\hat{t}_0(t_{min})$ as the final estimate of $t_0.$
\item
For the remaining simulation time ($t > t_{min}$), switch from the full system to the reduced model. The reduced model is evolved with the necessary parameter $t_0$ set to its estimated value $\hat{t}_0(t_{min}).$  
\end{enumerate}

This procedure can be extended to the computation of optimal estimates for $t_1,t_2,\ldots,$ i.e. when we evolve, in addition to $w_{0k}(t),$ the quantities $w_{1k}(t),w_{2k}(t),\ldots.$ Results for such higher order models will be presented elsewhere.


\section{Numerical Examples}\label{examples}

\subsection{A linear ODE with uncertain coefficient}\label{example1}
Consider the following {\it linear} ordinary equation with an uncertain coefficient
\begin{equation}\label{ex:ODE}
\begin{split}
\frac{du}{dt} &= -\kappa u,\\
u(0,\cdot )& = u^\circ,
\end{split}
\end{equation}
where $\kappa \sim U[0,1]$. This equation has the solution $u = u^\circ exp(-\kappa t)$. To represent the dependence of the solution of \eqref{ex:ODE} on $\kappa,$ we can expand it in a general polynomial chaos (gPC) expansion \cite{xiu2006}, say using Legendre polynomials. Let $u(t,\cdot)\approx \sum_{i=0}^M u_i(t) \phi_i(\xi)$, where $\xi \sim U[-1,1]$ and $\{\phi_i\}$ are normalized Legendre polynomials which are orthonormal with respect to the uniform distribution of $\xi$, i.e., $$\int_{-1}^1 \phi_i(\xi)\phi_j(\xi)\frac{1}{2}d\xi = \delta_{ij}.$$ We can write $\kappa$ as $\kappa =\frac{1}{2}\xi+\frac{1}{2}= \sum_{i=0}^1{k_i}\phi_i(\xi)$. We substitute this expansion in \eqref{ex:ODE} and obtain (through Galerkin projection) the (truncated) system up to order $M$
\begin{equation}\label{ex:ODE_g}
\begin{split}
\frac{du_r}{dt}& = -\sum_{i=0}^1\sum_{j=0}^M k_iu_je_{ijr}, \\
u_r(0) &= u_{0r}, \qquad r = 0,\dots, M,
\end{split}
\end{equation}
where $e_{ijk} = \int_{-1}^1 \phi_i(\xi)\phi_j(\xi)\phi_k(\xi)\frac{1}{2}d\xi$ and $u_{00} = u^\circ$, $u_{0r} = 0$ for $r = 1,\dots,M$ (for details, see e.g \cite{xiu2002}).

Because of the spectral decay in the gPC coefficients, it is natural to choose the coefficients $u_i$, $i = 0,\dots,\Lambda$ of the lower degree Legendre polynomials to be the resolved variables $\hat{u},$ and $u_i$, $i = \Lambda+1,\dots, M$ to be the unresolved variables $\tilde{u}$ respectively. 
To conform with the notation in Section \ref{mz_formalism}, we have $H = \{0,\dots,\Lambda\}$ and $G = \{\Lambda+1,\dots,M\}$.
We have chosen $M = 6$ for the full system and $\Lambda = 1$ for the reduced system. The solution of the full system is converged for $M=6$ and thus, we do not need to keep further terms in the expansion.  

The projection $P$ we have chosen is defined as $(Pf)(\hat{u}_0) = f(\hat{u}_0,\tilde{0}).$  Also, we define $Q=I-P.$ To be consistent with the notation in Section \ref{mz_formalism}, we have
$$ R_r({u}_0) = -\sum_{i=0}^1\sum_{j=0}^M k_iu_je_{ijr}. $$
\[
PLu_{0r} = -\sum_{i=0}^1\sum_{j=0}^{\Lambda}k_iu_{0j}e_{ijr}.
\]
\[
QLu_{0r} = -\sum_{i=0}^1\sum_{j=\Lambda+1}^Mk_iu_{0j}e_{ijr},
\]

In order to be able to compute the expressions for the memory terms we use a finite-rank projection $\mathbb{P}$ to approximate the projection $P$. To define the finite-rank projection we need to introduce a measure for the distribution of the coefficients. We consider the coefficients $u_{0r}$ to be i.i.d Gaussian random variables with mean at the values given initially (see \eqref{ex:ODE_g}) and a prescribed variance for $i = 0,\dots, M.$ In the case of the linear ODE, the variance was set to 0.01 for all the variables in the full system.  Also, $\omega$ is the joint probability measure with respect to these random variables. Then for a function $\varphi_j(u_0,t)$ of the initial conditions and time, the finite-rank projection reads
\begin{equation}\label{def:frank_proj}
(\mathbb{P}\varphi_j)(\hat{u}_0,t) = \sum_{\nu\in I}(\varphi_j(u_0,t),h^{\nu}(\hat{u}_0))h^{\nu}(\hat{u}_0),
\end{equation}
where $h^{\nu}(\hat{u}_0)$ are tensor product Hermite polynomials up to some order $p$, $\nu$ is the multi-index $\nu = (\nu_0,\dots,\nu_{\Lambda})$ with $|\nu| = \sum_{i=0}^\Lambda \nu_i$ and $I$ is the index set up to order $p$, i.e., $I = \{ \mu \big| |\mu|\leq p \}$. The order $p$ for the basis functions was set to 5 for a total of 21 basis functions. In formula \eqref{def:frank_proj} the inner product is defined as
\begin{equation}
(f,g) = \int fg d\omega . 
\end{equation}
For each $j\leq \Lambda$, the component $F_j(u_0,t)$ denotes the solution of the orthogonal dynamics
\begin{equation}\label{Orth_j}
\begin{split}
&\frac{\partial}{\partial t}F_j(u_0,t) ={Q}LF_j(u_0,t) = {L}F_j(u_0,t)-{P}LF_j(u_0,t),\\
& F_j(u_0,0) = {Q}Lu_{0j} = R_j(u_0)-{P}Lu_{0j}.
\end{split}
\end{equation} 
\eqref{Orth_j} is equivalent to the Dyson formula:
\begin{equation}\label{Orth_j_Dyson}
F_j(u_0,t) = e^{tL}F_j(u_0 ,0)-\int_{0}^t e^{(t-s)L}{P}LF_j(u_0,s)ds.
\end{equation}
Eq. \eqref{Orth_j_Dyson} is a Volterra integral equation for $F_j(u_0,t).$ To proceed, we replace the projection operator ${P}$ with the finite-rank projection operator $\mathbb{P}$ and find
\begin{equation}
K_j(\hat{u}_0,s) = PLF_j(u_0,s) \approx \mathbb{P}LF_j(u_0,s) =  \sum_{\nu \in I}a_j^{\nu}h^{\nu}(\hat{u}_0),
\end{equation}
where 
\[
a^{\nu}_j(s) = (LF_j(u_0,s),h^{\nu}(\hat{u}_0)).
\]
Consequently,
\[
e^{(t-s)L}\mathbb{P}LF_j(u_0,s) = \sum_{\nu\in I}a^{\nu}_j(s)h^{\nu}(\varphi(u,t-s)).
\]
We substitute $e^{(t-s)L}\mathbb{P}LF_j(u_0,s)$ for $e^{(t-s)L}PLF_j(u_0,s)$ in Eq. \eqref{Orth_j_Dyson}, multiply both sides by $L$ and take the inner product with $h^{\mu}(\hat{u}_0))$; the result is (dropping the approximation sign)

\begin{equation}
\begin{split}
&(LF_j(u_0,t),h^{\mu}(\hat{u}_0)) \\
=& (Le^{tL}F_j(u_0,0),h^{\mu}(\hat{u}_0))-\int_{0}^t \sum_{\nu\in I} a^{\nu}_j(s)(Le^{(t-s)L}h^{\nu}(\hat{u}_0),h^{\mu}(\hat{u}_0))ds. \label{Volterra_a1}
\end{split}
\end{equation}
Eq. \eqref{Volterra_a1} is a Volterra integral equation for the function $a^{\nu}_j(t)$, which can be rewritten as follows:
\begin{equation}\label{Volterra_a}
a^{\mu}_j(t) = f^{\mu}_j(t)-\int_{0}^t\sum_{\nu\in I}a^{\nu}_j(s)g^{\nu\mu}(t-s)ds,
\end{equation}
where 
\[
f^{\mu}_j(t) = (Le^{tL}F_j(u_0,0),h^{\mu}(\hat{u}_0)), \qquad g^{\nu\mu}(t)=(Le^{tL}h^{\nu}(\hat{u}_0),h^{\mu}(\hat{u}_0)).
\]
The functions $f^{\nu}_j(t)$, $g^{\mu\nu}(t)$ can be found by averaging over a collection of experiments or simulations, with initial conditions drawn from the initial distribution. In this example, we use a sparse grid quadrature rule for the multi-dimensional integrals \cite{xiu2006}.  

Finally, we perform one more projection to eliminate the noise term (see Section \ref{mz_formalism}) and the memory term becomes
\[
\int_{0}^t Pe^{(t-s)L}K_j(\hat{u}_0,s)ds.
\]
This can be approximated by
\[
\int_{0}^t \sum_{\nu\mu\in I}a^{\nu}_j(s)\gamma^{\nu\mu}(t-s)h^{\mu}(\hat{u}_0)ds,
\]
where 
\[
\gamma^{\nu\mu}(t) = (e^{tL}h^{\nu}(\hat{u}_0),h^{\mu}(\hat{u}_0)).
\]
After calculating $a_i^{\mu}$ and $\gamma^{\mu\nu}$ we obtain the following reduced system,
\begin{equation}\label{eq:redu_b}
\frac{d}{dt}\hat{u}(t) = \textrm{R}(\hat{u}(t))+\int_{0}^t A(s)\Gamma(t-s)h(\hat{u}_0)ds. \quad \hat{u}(0) = \hat{u}_0,
\end{equation}
here $A$ and $\Gamma$ are the matrix form of $a^{\mu}_i$ and $\gamma^{\mu\nu}$, $\hat{u}_0$ is the initial condition of resolved variables.

\begin{figure}[htbp]
\centering
\psfig{file = 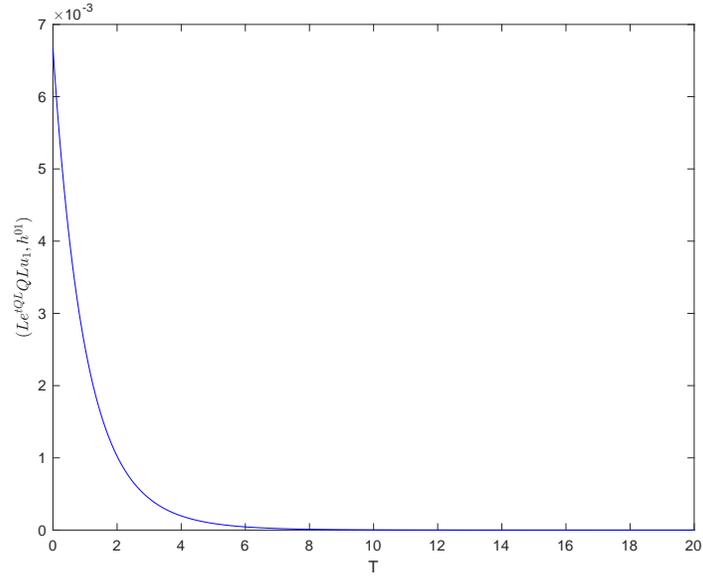, width = 11cm}
\caption{Evolution of the memory kernel $(Le^{tQL}QLu_1,h^{01})$ (see text for details).}
\label{fig:ODE_1_me}
\end{figure}

\begin{figure}[htbp]
\centerline{
\psfig{file = 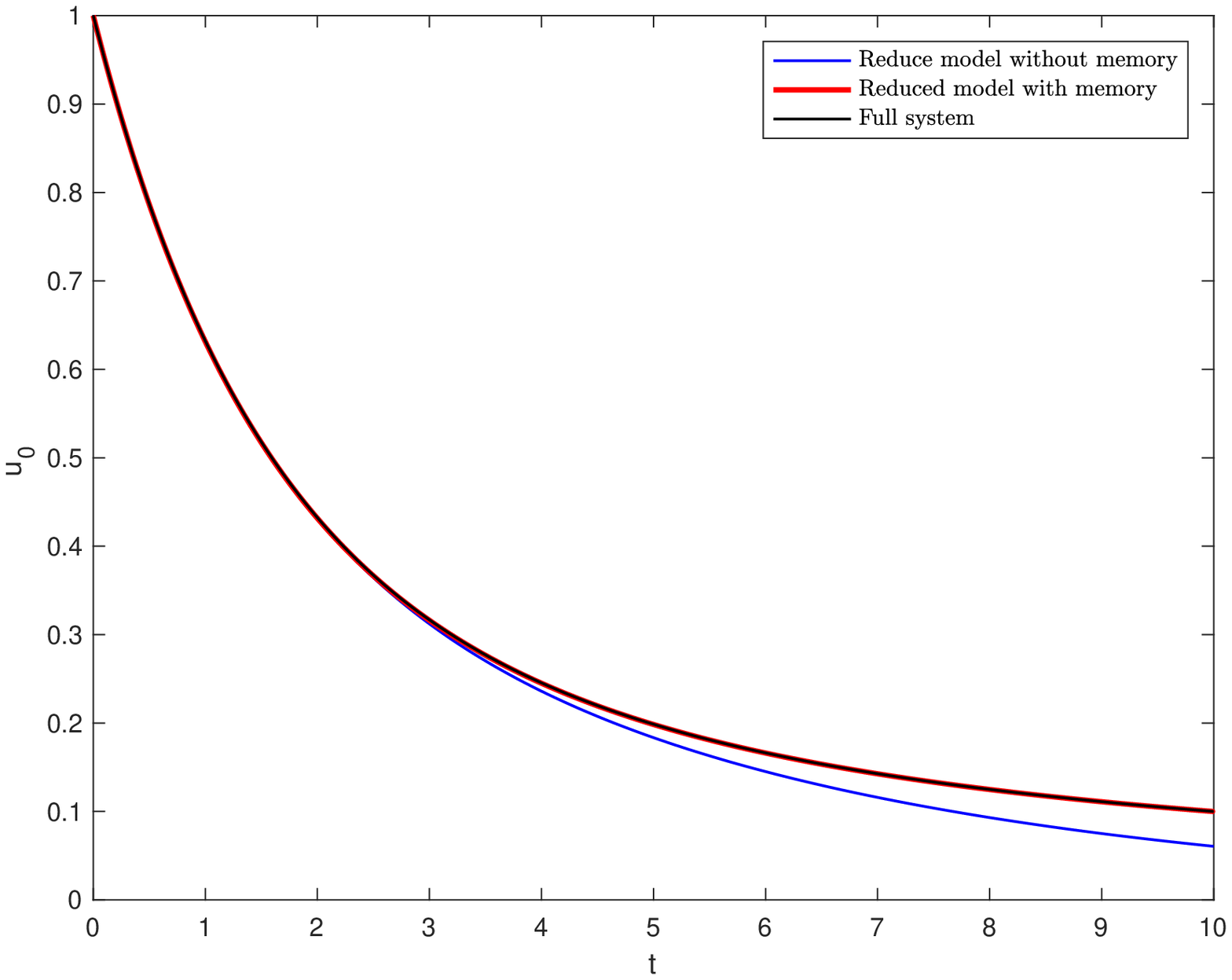, width = 7cm}
\psfig{file = 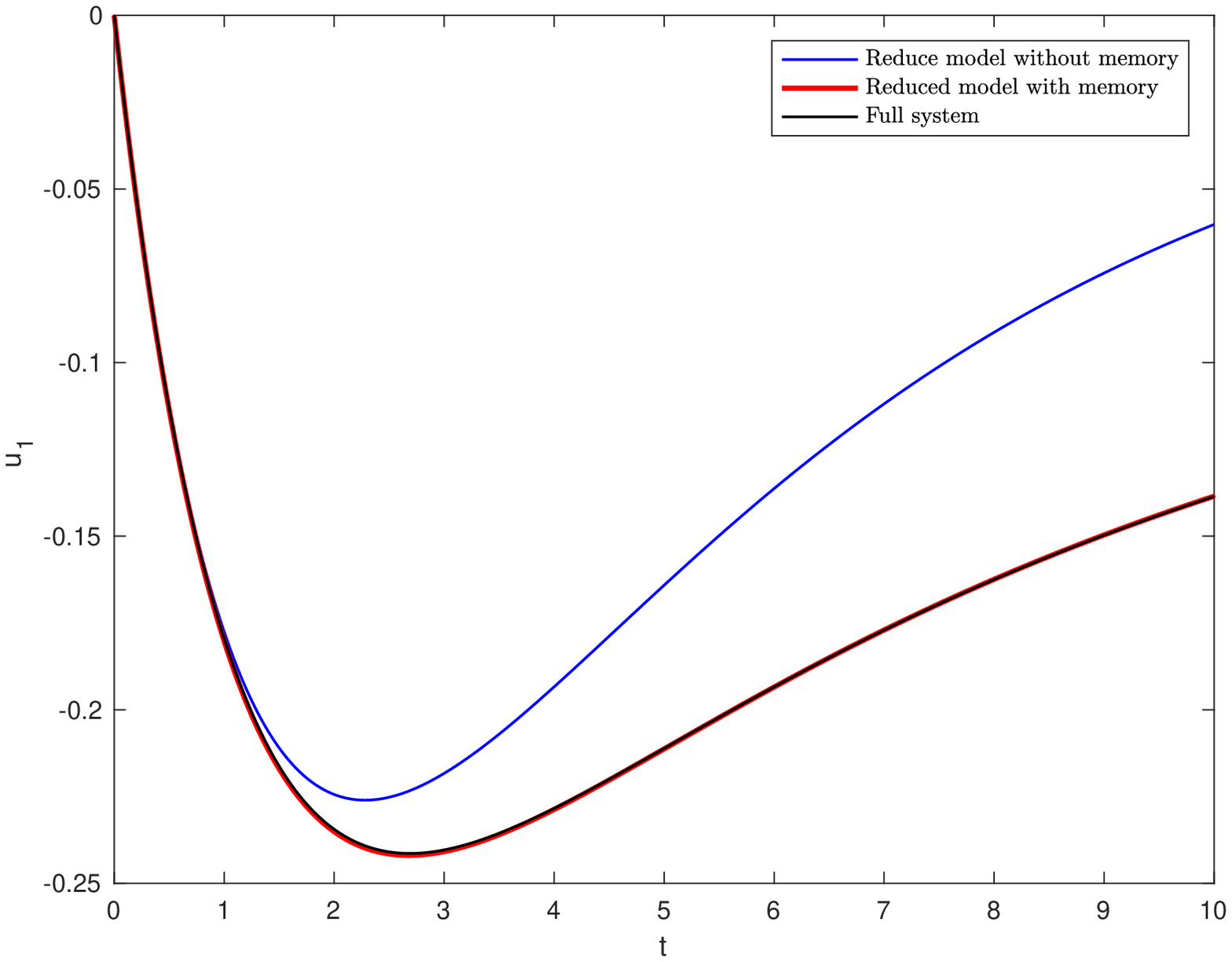, width = 7cm}
}
\caption{Evolution of the resolved variables $u_0,u_1$ predicted by the full model (black line), the (Markovian) reduced model without memory (blue line) and the (non-Markovian) reduced model with memory (red line).}
\label{fig:ODE_1_sol}
\end{figure}

\begin{figure}[htbp]
\centerline{
\psfig{file = 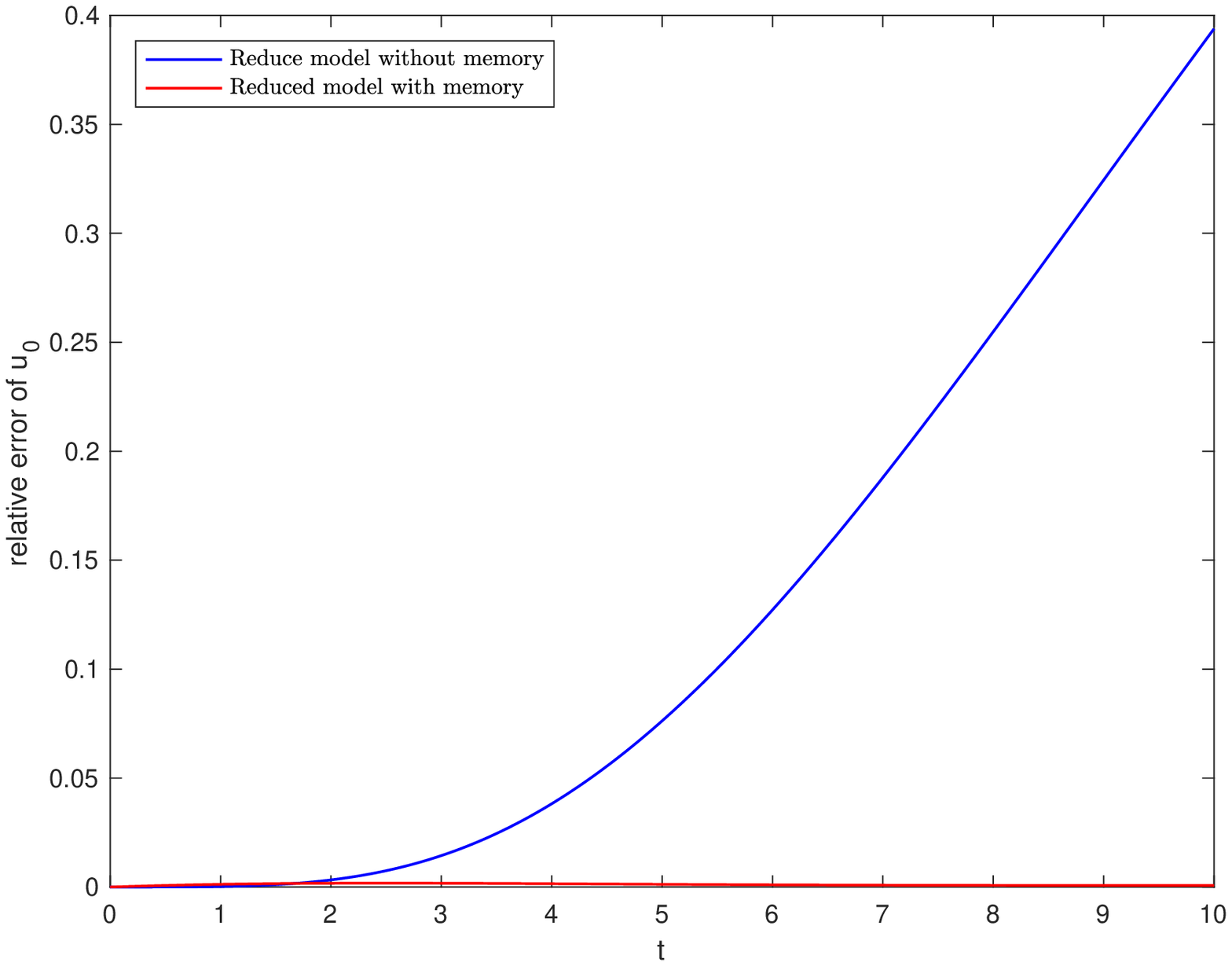, width = 7cm}
\psfig{file = 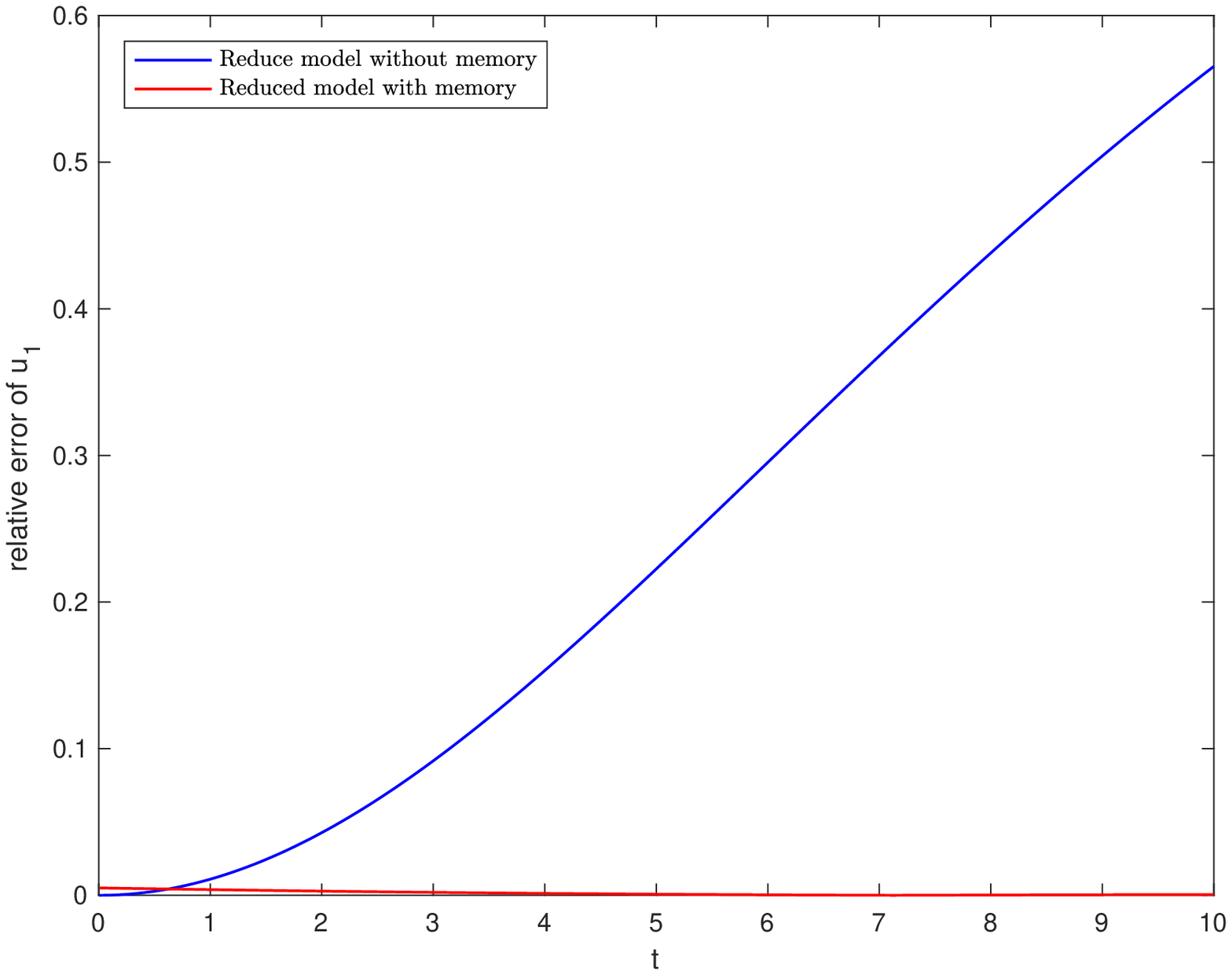, width = 7cm}
}
\caption{Relative error with respect to the true solution for the (Markovian) reduced model without memory (blue line) and the (non-Markovian) reduced model with memory (red line).}
\label{fig:ODE_1_er}
\end{figure}

Fig. \ref{fig:ODE_1_me} shows the evolution of the memory kernel $(Le^{tQL}QLu_1,h^{01})$ which is indicative of the behavior of the memory kernels. The basis function $h^{01}$ is the product of the zero order Hermite polynomial in the variable $u_0$ and the first order Hermite polynomial in the variable $u_1.$ We see that the memory kernel is rather slowly decaying which means that the resulting reduced order model will have a long memory. Fig. \ref{fig:ODE_1_sol} shows the solution for the resolved variables as predicted by the full system and two different reduced order models, the Markovian model which results from dropping the memory term in \eqref{eq:redu_b} and the non-Markovian reduced model given by \eqref{eq:redu_b}. It is obvious from Fig. \ref{fig:ODE_1_sol} that the Markovian model loses accuracy quickly. On the other hand, the non-Markovian model retains its accuracy for the length of the simulation interval. This difference in behavior is quantified in Fig. \ref{fig:ODE_1_er} where we see that for both resolved variables the relative error of the Markovian model becomes greater than $50\%$ by the end of the simulation interval. On the other hand, the error of the non-Markovian model remains less than $1\%$ for the whole simulation interval.


\subsection{Nonlinearly damped and randomly forced particle}\label{example2}
Consider the following equation describing a particle moving in a double well potential and driven by a force term (see \cite{xiu2006})
\begin{equation}\label{eq:Damped_system}
\begin{split}
\frac{du}{dt} & = u-u^3+f(t,\xi),\\
u(0) & = u^\circ,\\
\end{split}
\end{equation}
where $f = \sin(t+t_0)\xi$ and $\xi \sim U[-1,1]$. We use order $M=6$ polynomials in $\xi$ to approximate the full system solution up to time $10$. 
We want to construct a reduced model for the first 2 coefficients of the polynomial expansion ($\Lambda=1$). As before, we let $u(t,\xi)\approx \sum_{i=0}^{M} u_i(t)\phi_i(\xi)$ and we obtain through Galerkin projection the system
\begin{equation}
\begin{split}
\frac{du_i}{dt} &= u_i-\sum_{j,k,m = 0}^{M}u_j u_k u_m e_{jkmi}+f_i,\\
u_i(0) &= u_{0i}, \quad \textrm{ for } i = 0,\dots,M.
\end{split}
\end{equation} 
where $e_{jkmi} = \int \phi_j(\xi) \phi_k(\xi) \phi_m(\xi) \phi_i(\xi) \frac{1}{2}d\xi$ and $u_{0i} = \{\begin{array}{ll} u^\circ, & i=0;\\ 0,& \textrm{otherwise}. \end{array}$
Let $R$ be the vector with $R_i = u_i-\sum_{j,k,m = 0}^{M}u_j u_k u_m e_{jkmi}+f_i $. In order to apply the MZ formalism we need an autonomous system of equations to begin with. For this purpose, we introduce an auxiliary time-variable $\tau,$ such that $\tau = t$ and $\frac{d\tau}{dt} = 1.$ The projection operator $\mathbb{P}$ projects onto the function space of the first two coefficients and $\tau$. Again, we use the finite-rank projection onto the function space expanded by Hermite polynomials up to order $3$ to represent the orthogonal dynamics (total of 10 functions) and solve the Volterra equation for the memory kernels as we did for the linear ODE example. The variance for the Gaussian variables used to define the inner product for the finite-rank projection was set to $10^{-2i-2}$ for the coefficient $u_i$ with $i=0,1,\ldots,6.$ The reason we used a decreasing sequence of variances as we go up in the order of Legendre polynomials is to stabilize the behavior of the reduced model.

\begin{figure}[htbp]
\centerline{
\psfig{file = 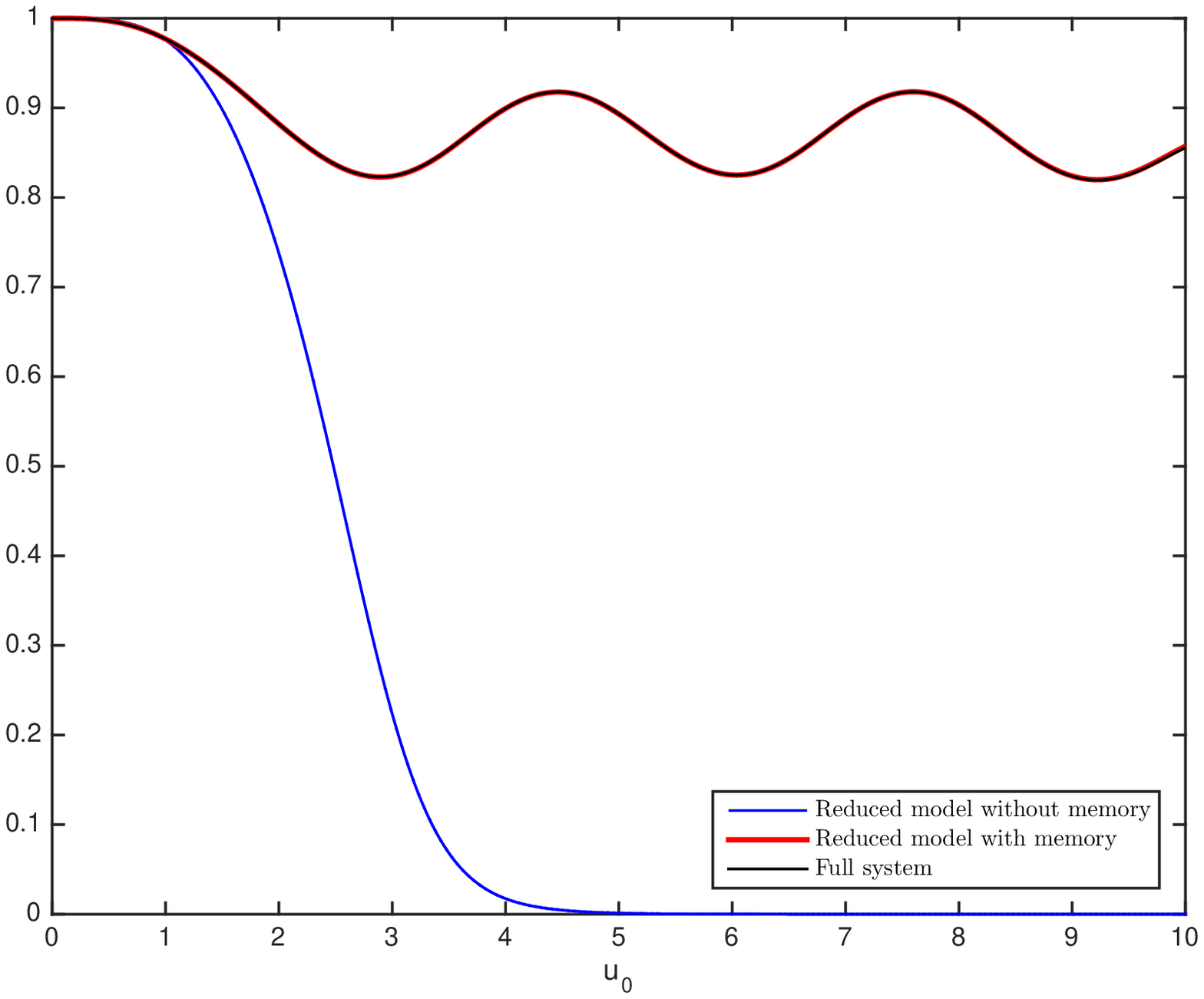, width = 7cm}
\psfig{file = 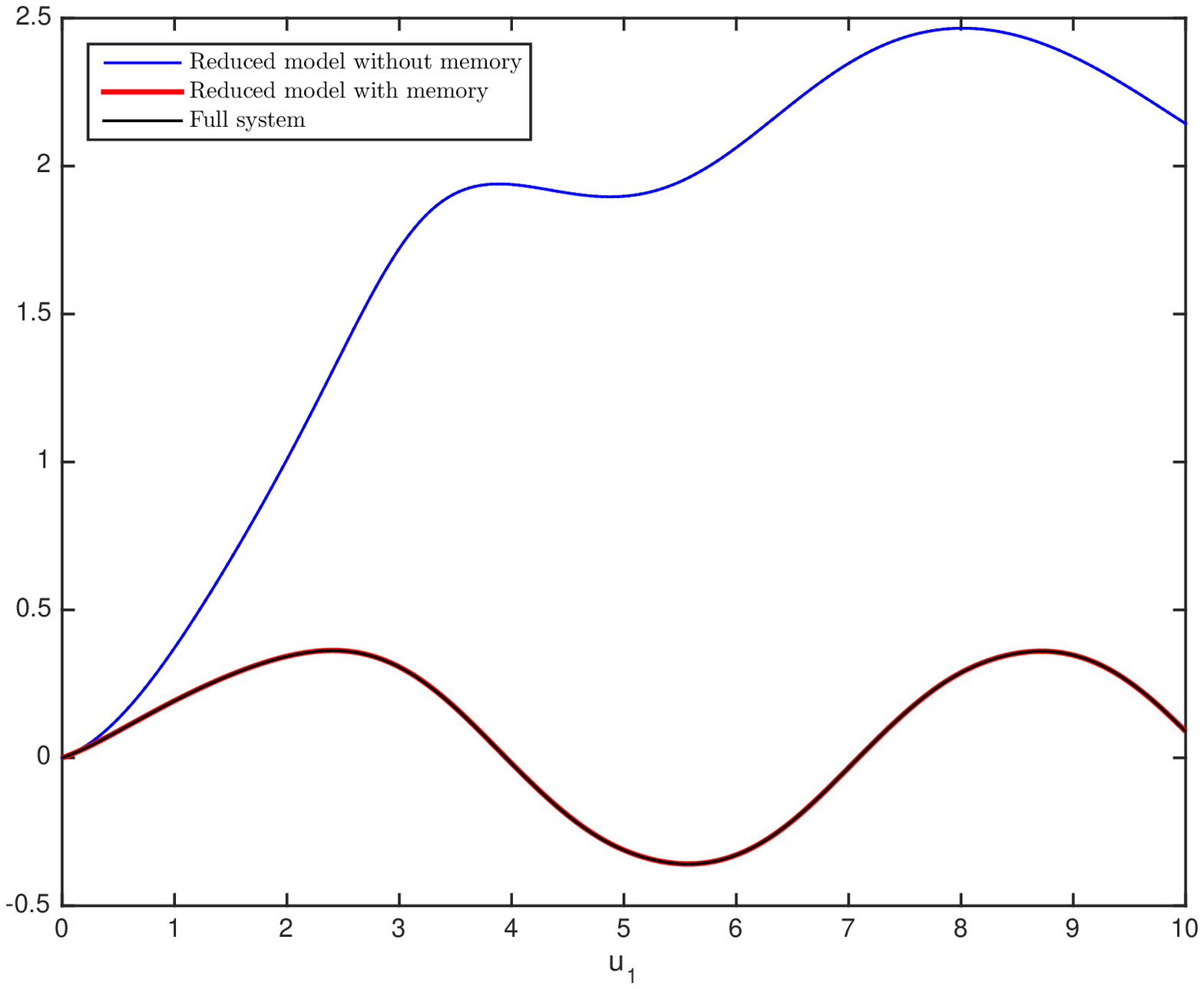, width = 7cm}
}
\caption{Evolution of the resolved variables $u_0,u_1$ predicted by the full model (black line), the (Markovian) reduced model without memory (blue line) and the (non-Markovian) reduced model with memory (red line).}
\label{fig:damped_T10}
\end{figure}

\begin{figure}[htbp]
\centerline{
\psfig{file = 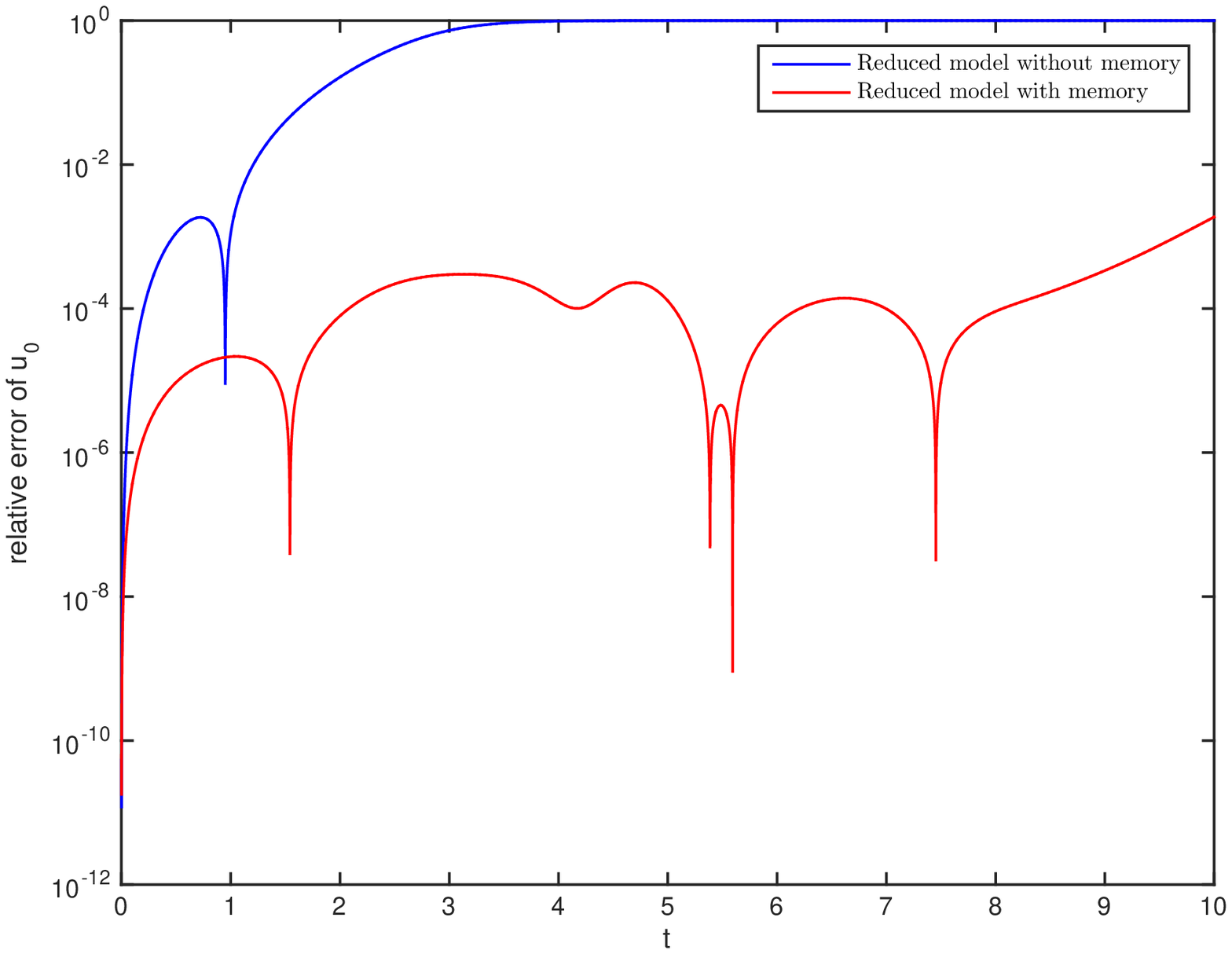, width = 7cm}
\psfig{file = 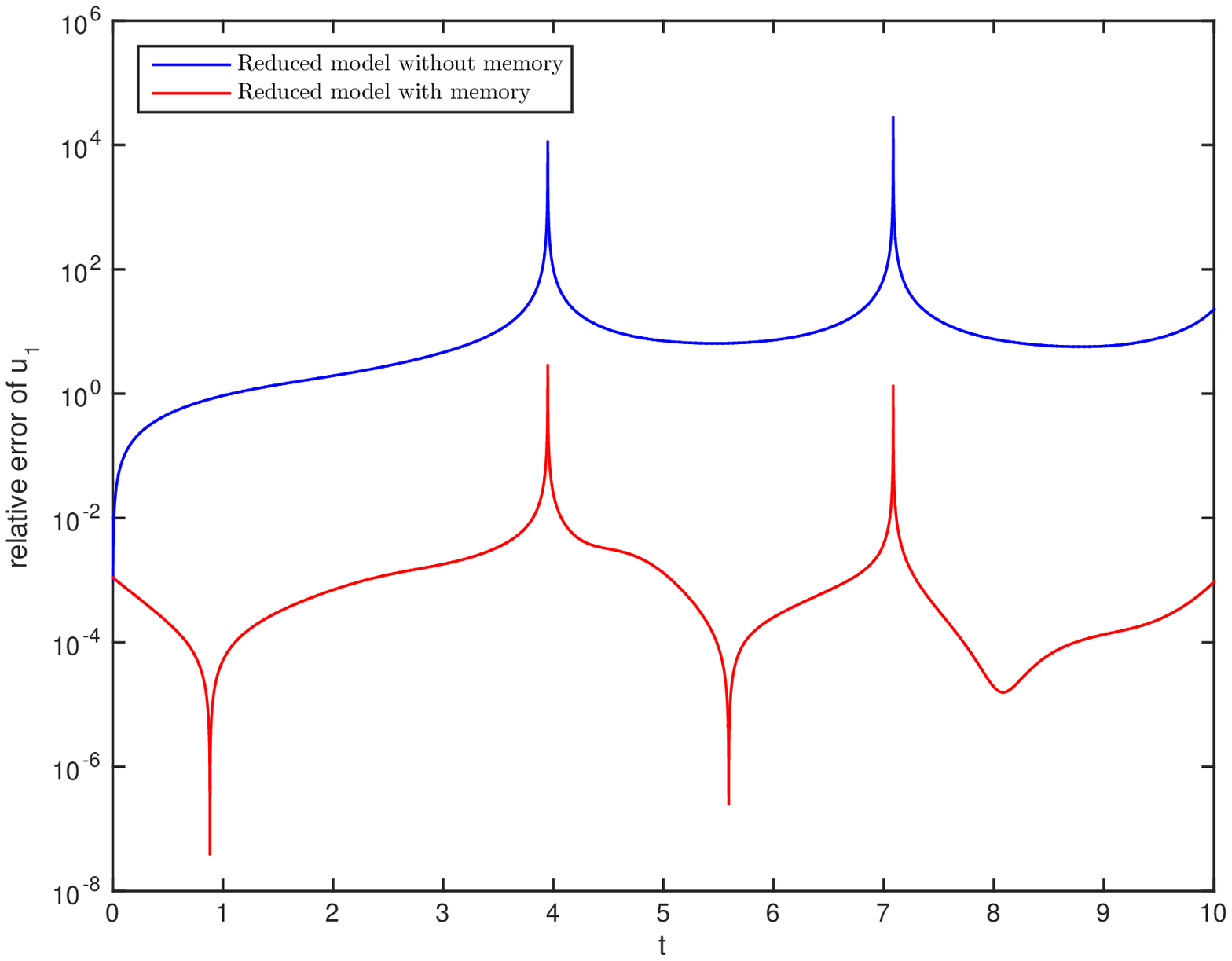, width = 7cm}
}
\caption{Logarithimic scale relative error for $u_0,u_1$ with respect to the true solution for the (Markovian) reduced model without memory (blue line) and the (non-Markovian) reduced model with memory (red line).}
\label{fig:damped_T10_log_er}
\end{figure}

As can be seen from Fig. \ref{fig:damped_T10_log_er}, the difference between the (memoryless) Markovian and non-Markovian reduced models is even more pronounced than in the case of the linear ODE. The inclusion of the memory term is indeed crucial for maintaining the accuracy of the reduced model for long times. For the case of the resolved variable $u_1,$ the relative error spikes at a couple of points even for the otherwise very accurate non-Markovian reduced model. As can be seen from Fig. \ref{fig:damped_T10}, this is because the exact value of $u_1$ becomes zero at these points so that the relative error becomes very large even for an accurate approximation. However, the significant improvement in accuracy with the inclusion of the memory term is evident in Fig. \ref{fig:damped_T10_log_er} which plots the error in a logarithmic scale.   

\subsection{Viscous 1D Burgers with uncertain initial conditions}\label{example3}

In this section we show how the above MZ formulation can be used for uncertainty quantification for the one-dimensional Burgers equation with uncertain initial condition. As is explained at the end of this section, the calculation of the MZ memory term cannot proceed as for the last two examples. The reason is that it is prohibitively expensive due to the number of basis functions needed. Thus, we will apply the alternative construction that was presented in Section \ref{reformulation}.

The equation is given by 
\begin{equation}\label{burgersequation}
u_t+u u_x = \nu u_{xx},
\end{equation}
where $\nu > 0.$ Equation (\ref{burgersequation}) should be supplemented with an initial condition $u(x,0)=u_0(x)$ and boundary conditions. We solve (\ref{burgersequation}) in the interval $[0,2\pi]$ with periodic boundary conditions. This allows us to expand the solution in Fourier series
$$u_{N}(x,t )=\underset{k \in F}{\sum} u_k(t) e^{ikx},$$
where $F=[-\frac{N}{2},\frac{N}{2}-1].$ The equation of motion for the Fourier mode $u_k$ becomes
\begin{equation}
\label{burgersode}
 \frac{d u_k}{dt}=- \frac{ik}{2} \underset{p, q \in F}{\underset{p+q=k  }{ \sum}} u_{p} u_{q}  -\nu k^2 u_k.
\end{equation}
We assume that the initial condition $u_0(x)$ is uncertain (random) and can be expanded as $u_0(x,\xi)= (\alpha_0+ \alpha_1 \xi) v_0(x)$ where $\xi$ is uniformly distributed in $[-1,1]$ and $v_0(x)$ a given function. In the numerical experiments we have taken $\alpha_0 =\alpha_1=1$ and $v_0(x)=\sin x.$ Thus, the initial condition varies ``uniformly" between the functions 0 and $2 \sin x.$  

To proceed we expand the solution $u_k(t,\xi)$ for $k \in F$ in a polynomial chaos expansion using the standard Legendre polynomials which are orthogonal in the interval $[-1,1].$ In particular, we have that $$\int_{-1}^1\phi_i (\xi) \phi_j(\xi) \frac{1}{2}d \xi=\frac{1}{2i+1} \delta_{ij},$$ where $\phi_i(\xi)$ is the standard Legendre polynomial of order $i.$ For each wavenumber $k$ we expand the solution $u_k(t,\xi)$ of \eqref{burgersodemz} in Legendre polynomials and keep the first $M$ polynomials  
\begin{equation}\label{ode_expansion}
u_k(t,\xi)\approx \sum_{i=0}^{M-1} u_{ki}(t) \phi_i(\xi), \; \; \text{where} \; \; \xi \sim U[-1,1].
\end{equation}
Similarly, the initial condition can be written as $u_0(x,\xi) = \sin x \sum_{i=0}^1\alpha_i \phi_i(\xi)$ since $\phi_0 (\xi)=1$ and $\phi_1 (\xi)=\xi.$
Substitution of \eqref{ode_expansion} in \eqref{burgersode} and use of the orthogonality property of the Legendre polynomials gives 
\begin{equation}\label{burgersodemz_system}
\frac{du_{kr}(t)}{dt}=- \frac{ik}{2} \sum_{l=0}^{M-1} \sum_{m=0}^{M-1} \underset{p, q \in F}{\underset{p+q=k  }{ \sum}} u_{pl} u_{qm}  c_{lmr} - \nu k^2 u_{kr}
\end{equation}
for $k \in F$ and $r=0,\ldots,M-1.$ Also $$c_{lmr}=\frac{E[ \phi_l (\xi) \phi_m(\xi) \phi_r(\xi) ] }{E[\phi^2_r(\xi) ]},$$
where the expectation $E[\cdot]$ is taken with respect to the uniform density on $[-1,1].$ The Legendre polynomial triple product integral defines a tensor which has the following sparsity pattern: $E[ \phi_l (\xi) \phi_m(\xi) \phi_r(\xi) ]=0,$ if $ l+m < r$ or $l+r < m$ or $m+r < l$ or $l+m+r= \text{odd}$ \cite{gupta}. Due to this sparsity pattern, for a given value of $M$ only about $1/4$ of the $M^3$ tensor entries are different from zero.

Before we proceed we have to comment on the cost of applying the MZ formalism to construct a reduced model. We have set the viscosity coefficient to $\nu = 0.03.$ The solution of the full system was computed with $N=196$ Fourier modes ($F=[-98,97]$) and the first 7 Legendre polynomials ($M=7$). The first 7 Legendre polynomials were enough to obtain converged statistics for the full system. We want to construct reduced models for the evolution of the coefficients of the first 2 Legendre polynomials i.e., $u_{k0},u_{k1}$ for $k \in F.$ If we want to apply the MZ formalism in the way we did for the previous two examples (employing a finite-rank projection etc.) we would need to construct a basis in $2 \times 98$ dimensions (exploiting the fact that the solution of the Burgers equation is real-valued). Any attempt to use basis functions up to a high order is infeasible for such a high-dimensional situation. We have attempted to use only low order basis functions but they are not enough to guarantee accuracy of the reduced model. Thus, we turn to the reformulated reduced model that was presented in Section \ref{reformulation}.

\subsubsection{Reformulated MZ reduced model}\label{mz_ode_example}

To conform with 
the Mori-Zwanzig formalism we set 
$$R_{kr}(u)=- \frac{ik}{2} \sum_{l=0}^{M-1} \sum_{m=0}^{M-1} \underset{p, q \in F}{\underset{p+q=k  }{ \sum}} u_{pl} u_{qm}  c_{lmr} - \nu k^2  u_{kr} ,$$
where $u=\{u_{kr}\}$ for $k \in F$ and $r=0,\ldots,M-1.$ Thus, we have
\begin{equation}
\label{burgersodemz}
\frac{d u_{kr}}{dt}=R_{kr}(u) 
\end{equation}
for $k \in F$ and $r=0,\ldots,M-1.$  
We proceed by dividing the variables in resolved and unresolved. In particular, we consider as resolved the variables $\hat{u}=\{u_{kr}\}$ for $k \in F$ and $r=0,\ldots,\Lambda-1,$ where $\Lambda < M.$ Similarly, the unresolved variables are $\tilde{u}=\{u_{kr}\}$ for $k \in F$ and $r=\Lambda,\ldots,M-1.$ In the notation of Section \ref{mz_formalism} we have $H= F \cup (0,\ldots,\Lambda-1)$ and $G= F\cup (\Lambda,\ldots,M-1).$ In other words, we resolve, for all the Fourier modes, only the first $\Lambda$ of the Legendre expansion coefficients and we shall construct a reduced model for them.

The system (\ref{burgersodemz}) is supplemented by the initial 
condition $u_0=(\hat{u}_0,\tilde{u}_0).$ We focus on initial conditions where 
the unresolved Fourier modes are set to zero, i.e. $u_0=(\hat{u}_0,0).$ We also define $L$ by 
$$L=\sum_{k \in F}\sum_{r=0}^{M-1} R_{kr}(u_0) \frac{\partial}{\partial u_{0kr}}.$$ 
To construct a MZ reduced model we need to define a projection operator $P.$ For a function $h(u_0)$ of all the 
variables, the projection operator we will use is defined by $P(h(u))=P(h(\hat{u}_0,\tilde{u}_0))=h(\hat{u}_0,0),$ i.e. 
it replaces the value of the unresolved variables $\tilde{u}_0$ in any function $h(u_0)$ by zero. Note that this choice of projection is consistent with the initial conditions we have chosen. Also, we define the Markovian term 
$$ PLu_{0k}=PR_k(u_0)=- \frac{ik}{2} \sum_{l=0}^{\Lambda-1} \sum_{m=0}^{\Lambda-1} \underset{p, q \in F}{\underset{p+q=k  }{ \sum}} u_{0pl} u_{0qm}  c_{lmr} - \nu k^2  u_{0kr}.$$ 
The Markovian term has the same functional form as the RHS of the full system but is restricted to a sum over only the first $\Lambda$ Legendre expansion coefficients  for each Fourier mode. 

For the the term $PLQLu_{0kr}$ we find

\begin{equation}\label{burgersmemory1}
PLQLu_{0kr}=2\times \biggl [   - \frac{ik}{2}   \sum_{l=\Lambda}^{M-1} \sum_{m=0}^{\Lambda-1} \underset{p, q \in F}{\underset{p+q=k  }{ \sum}} PLu_{0pl} u_{0qm}  c_{lmr} \biggr ] .
\end{equation}

Finally, to implement any method to solve equation \eqref{newton1} for the estimation of $t_0$ we need to specify the RHS of the equation \eqref{newton1}. This requires the evaluation of the expression $Pe^{tL}QLu_{0kr}.$ For the case of the viscous Burgers equation, we find
\begin{gather}\label{newton4}
Pe^{tL}QLu_{0kr}=2(- \frac{ik}{2}) \sum_{l=\Lambda}^{M-1} \sum_{m=0}^{\Lambda-1} \underset{p, q \in F}{\underset{p+q=k  }{ \sum}} u_{pl} u_{qm}  c_{lmr} \\
- \frac{ik}{2} \sum_{l=\Lambda}^{M-1} \sum_{m=\Lambda}^{M-1} \underset{p, q \in F}{\underset{p+q=k  }{ \sum}} u_{pl} u_{qm}  c_{lmr}. \notag
\end{gather}  
Note that since we restrict attention to initial conditions for which the unresolved variables are zero and the projection sets the unresolved variables to zero, the quantity $Pe^{tL}QLu_{0kr}$ can be computed through the evolution of the full system \eqref{burgersodemz}.

The full system was solved with the modified Euler method with $\delta t = 0.001.$ The reduced model uses $N=196$ Fourier modes but only the first two Legendre polynomials, so $\Lambda=2.$ It was solved using the modified Euler method with $\delta t = 0.001.$ The parameter $t_0$ needed for the evolution of the memory term was found to be 0.3783 through the procedure described in Section \ref{mz_optimal}.

\begin{figure}
\centering
\epsfig{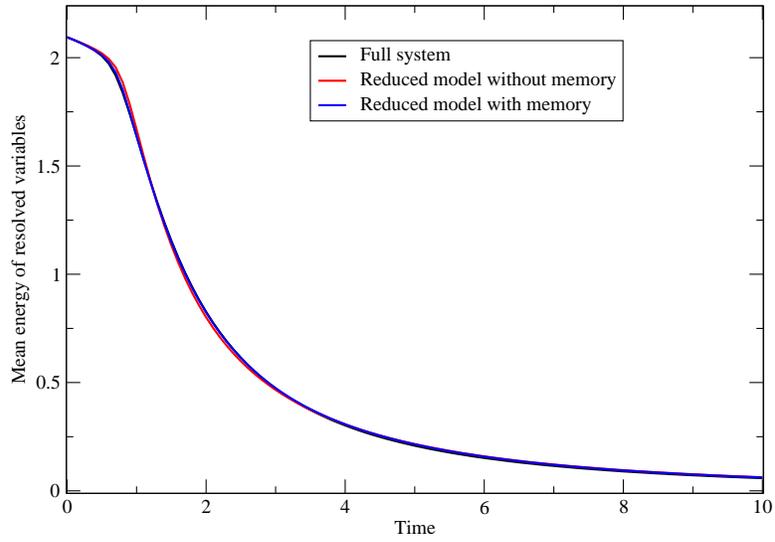}
\caption{Evolution of the mean of the energy of the solution using only the first two Legendre polynomials.}
\label{plot_initial_energy_mean}
\end{figure}

\begin{figure}
\centering
\epsfig{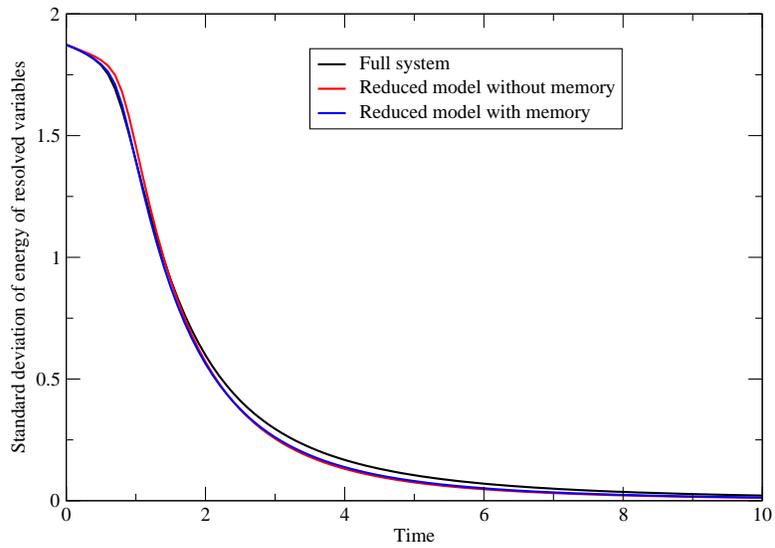}
\caption{Evolution of the standard deviation of the energy of the solution using only the first two Legendre polynomials.}
\label{plot_initial_energy_stdev}
\end{figure}

Figure \ref{plot_initial_energy_mean}  shows the evolution of the mean energy of the solution 
$$\mathbb{E}[E(t)]=\frac{1}{2}  \sum_{k \in F} \sum_{r=0}^1 2\pi |u_{kr}(t)|^2 \frac{1}{2r+1}$$
as computed from the full system (with $M=7$ Legendre polynomials), the MZ reduced model with $\Lambda=2$ {\it without} memory (keeping only the Markovian term) and the MZ reduced model with $\Lambda=2$ {\it with} memory. Figure \ref{plot_initial_energy_stdev} shows the evolution of the standard deviation of the energy of the solution. The variance of the energy is given by
$$Var[E(t)]=\frac{1}{4}  \sum_{k_1, k_2 \in F}   \sum_{r_1,\ldots,r_4=0}^1 (2\pi)^2 u_{k_1r_1} u_{k_1r_2}^*u_{k_2r_3} u_{k_2r_4}^*d_{r_1r_2r_3r_4}-\{ \mathbb{E}[E(t)]\}^2,$$
where
$$d_{r_1r_2r_3r_4}=\int_{-1}^1\phi_{r_1}(\xi)\phi_{r_2}(\xi)\phi_{r_3}(\xi)\phi_{r_4}(\xi) \frac{1}{2} d\xi.$$
The reduced model performs equally well with or without memory. Of course, the reduced model with memory is slower than the reduced model without memory. However, the reduced model with memory is still about 4 times faster than the the full system.

\begin{figure}
\centering
\epsfig{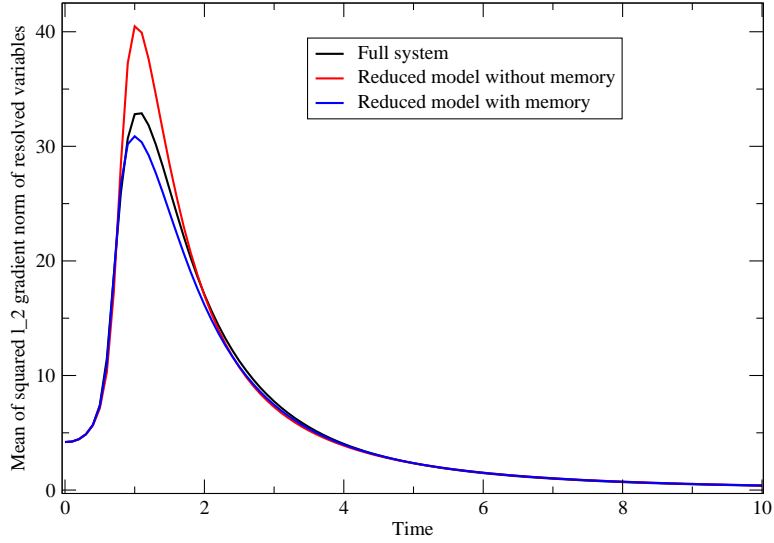}
\caption{Evolution of the mean of the squared $l_2$ norm of the gradient of the solution calculated using only the first two Legendre polynomials.}
\label{plot_initial_gradient_mean}
\end{figure}

\begin{figure}
\centering
\epsfig{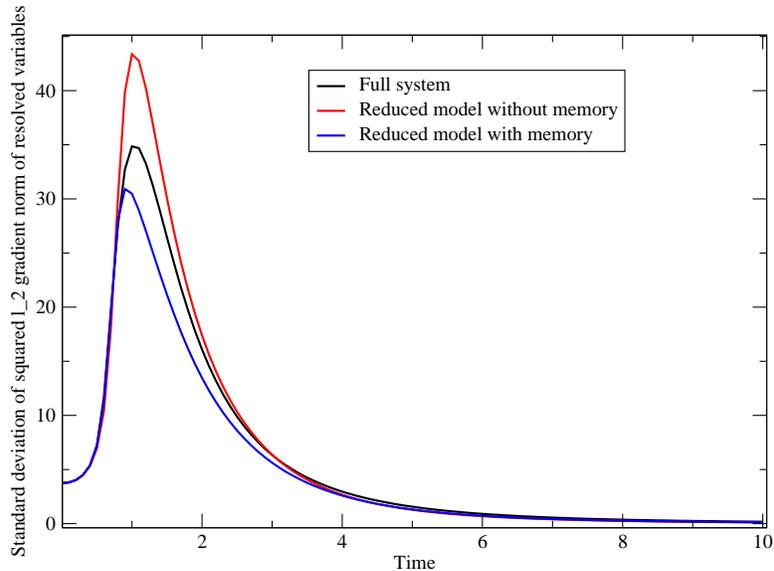}
\caption{Evolution of the standard deviation of the squared $l_2$ norm of the gradient of the solution calculated using only the first two Legendre polynomials.}
\label{plot_initial_gradient_stdev}
\end{figure}

Figure \ref{plot_initial_gradient_mean} shows the evolution of the mean squared $l_2$ norm of the gradient of the solution 
$$\mathbb{E}[G(t)]= \sum_{k \in F} \sum_{r=0}^1 2\pi k^2 |u_{kr}(t)|^2 \frac{1}{2r+1}$$
as computed from the full system (with $M=7$ Legendre polynomials), the MZ reduced model with $\Lambda=1$ {\it without} memory (keeping only the Markovian term) and the MZ reduced model with $\Lambda=1$ {\it with} memory. Figure \ref{plot_initial_gradient_stdev}  shows the evolution of the standard deviation. The variance is given by $$Var[G(t)]= \sum_{k_1, k_2 \in F}   \sum_{r_1,\ldots,r_4=0}^1 (2\pi)^2 k_1^2 k_2^2u_{k_1r_1} u_{k_1r_2}^*u_{k_2r_3} u_{k_2r_4}^*d_{r_1r_2r_3r_4}-\{ \mathbb{E}[G(t)]\}^2.$$
The large values of the standard deviation of the mean squared $l_2$ norm of the gradient are justified by the uncertainty in the initial condition. Recall that we have chosen an initial condition which can vary ``uniformly" between the functions 0 and $2\sin x.$ As a result, the standard deviation is large because it has to account for a wide range of possible initial conditions. 

It is evident from the figures that the inclusion of the memory term improves the performance of the reduced model. Also, it is evident that there is room for improvement of the reduced model {\it with} memory. In particular, more terms are needed in the reformulated MZ model to approximate better the memory. 

Recall that the solution of Burgers equation is a contraction \cite{lax}. Eventually, the complete description of the uncertainty caused by the uncertainty in the initial condition requires only a few polynomial chaos expansion coefficients. This happens at a time scale that is dictated by the magnitude of the viscosity coefficient. That is why for long times the reduced model with and without memory have comparable behavior to that of the full system. However, for short times, the inclusion of the memory term does make a difference because information from the higher chaos expansion coefficients is needed. The higher chaos expansion coefficients will have a more prolonged contribution for systems that possess unstable modes. In such cases, the inclusion of the memory term becomes imperative for short as well long times. Results for such cases will be presented elsewhere.


\section{Discussion and future work}\label{discussion}

We have examined the application of the Mori-Zwanzig formalism to the problem of constructing reduced models for uncertainty quantification. In particular, we have constructed reduced models for subsets of the polynomial chaos expansion coefficients needed to describe fully the uncertainty. We have examined cases of parametric or initial condition uncertainty. The main conclusion from the current work is that while the MZ formalism can be applied for the construction of reduced models, the task of constructing an efficient (or even feasible) reduced model can be involved. For cases where the straightforward application of the MZ formalism is not possible, we have offered an alternative construction. The implementation of this alternative construction is reminiscent of renormalization constructions used to describe the evolution of complex solutions of PDEs \cite{s11}.

The current work opens several directions for future work. First, we should investigate whether there is a more economical way of choosing the basis functions for cases when the basis functions have many arguments (as was the case for the Burgers example). This is important because the calculation of the memory kernels through the finite-rank projection is well defined and the solution of the corresponding Volterra equations can be performed with high accuracy. A related question is whether there is sparsity in the coefficients of the basis functions. It is plausible that even though in principle the number of basis functions to reach a specific order may be very large, many of them may not contribute to the representation. A related approach would be the use of machine learning algorithms to obtain a more efficient representation of the memory term. Finally, a related issue to be investigated is how to ensure the stability of the reduced model when the finite-rank projection is employed. For example, for the nonlinearly damped and forced particle case, we had to assign smaller variances for the higher coefficients to stabilize the reduced model. This procedure needs to be investigated and, if possible, automated.

A second interesting research direction has to do with the representation of the memory term when the finite-rank projection is not possible due to a prohibitively large number of basis functions. We have explored here an expansion of the memory term that involves, in essence, a Taylor expansion of the orthogonal dynamics operator. Such an expansion seems more plausible when the timescale of the orthogonal dynamics is {\it slower} than that of the resolved variables. However, there is an alternative way of performing the expansion of the memory term that is more suited to the case when the orthogonal dynamics is {\it faster} than the resolved variables. Such an expansion leads to a Taylor expansion of the whole memory term, not just the orthogonal dynamics operator. If the memory kernel becomes insignificant after a time interval $t_0,$ then one can use the {\it full} system up to time $t_0,$ estimate the Taylor expansion of the whole memory term around time $t_0$ and then switch to the {\it reduced} model with the memory given by the Taylor expansion. We will investigate this alternative memory representation and report the results elsewhere.

\section{Acknowledgements}
The authors would like to thank D. Barajas-Solano, H. Lei and A. Tartakovsky for useful discussions and comments. This research at Pacific Northwest National Laboratory (PNNL) was partially supported by the U.S. Department of Energy (DOE) Office of Advanced Scientific Computing Research (ASCR) Collaboratory on Mathematics for Mesoscopic Modeling of Materials (CM4), under Award Number DE-SC0009280 and partially by the U.S. DOE ASCR project ``Uncertainty Quantification For Complex Systems Described by Stochastic Partial Differential Equations". PNNL is operated by Battelle for the DOE under Contract DE-AC05-76RL01830.

\end{document}